\documentclass[a4paper,11pt]{amsart}

\usepackage[centertags]{amsmath}
\usepackage{amsfonts}
\usepackage{amssymb}
\usepackage{amsthm}
\usepackage{newlfont}
\usepackage{euscript}

\newtheorem{thm}{Theorem}[section]
\newtheorem{lem}{Lemma}[section]

\newtheorem{prop}{Proposition}[section]
\theoremstyle{remark}
\newtheorem{zau}{Remark}[section]
\theoremstyle{definition}
\newtheorem{ozn}{Definition}[section]
\newtheorem{ex}{Example}[section]
 \numberwithin{equation}{section}

\begin{document}
\newcommand{\No}{N}
\renewcommand{\phi}{\varphi}
\newcommand{\1}{1\!\!{\mathrm I}}
\renewcommand{\Re}{{\Bbb R}}
\newcommand{\eps}{\varepsilon}
\newcommand{\kap}{\varkappa}
\newcommand{\vO}{\varOmega}
\newcommand{\bu}{\mathbf{u}}
\newcommand{\bp}{\mathbf{p}}
\newcommand{\bx}{{x_*}}
\newcommand{\bc}{\mathbf{c}}
\newcommand{\bt}{t_*}
\newcommand{\bq}{\mathbf{q}}
\newcommand{\bR}{\mathbf{R}}
\newcommand{\bY}{\mathbf{Y}}
\newcommand{\ax}{\Re^+}
\newcommand{\prt}{\partial}
\newcommand{\Es}{\mathsf{E}}
\newcommand{\Cs}{\mathsf{C}}
\newcommand{\Ps}{\mathsf{P}}
\newcommand{\Xs}{\mathsf{X}}
\newcommand{\Ys}{\mathsf{Y}}
\newcommand{\Ff}{{\EuScript F}}
\newcommand{\Xf}{\mathfrak{X}}
\newcommand{\Yf}{{\EuScript Y}}
\newcommand{\Ef}{{\EuScript E}}
\newcommand{\Bf}{\mathfrak{B}}
\newcommand{\Tf}{{\EuScript T}}
\newcommand{\Cf}{{\EuScript C}}
\newcommand{\Zf}{{\EuScript Z}}
\newcommand{\Wf}{{\EuScript W}}
\newcommand{\Df}{{\EuScript D}}
\newcommand{\Nf}{{\EuScript N}}
\newcommand{\Gf}{{\EuScript G}}
\newcommand{\Qf}{{\EuScript Q}}
\newcommand{\Af}{\mathfrak{A}}
\newcommand{\Lf}{{\EuScript L}}
\newcommand{\Jf}{{\EuScript J}}
\newcommand{\Uf}{{\EuScript U}}
\newcommand{\Sf}{{\EuScript S}}
\newcommand{\Pf}{{\EuScript P}}
\newcommand{\pf}{\Pi_{fin}}
\newcommand{\Kb}{\mathbf{K}}
\newcommand{\ZZ}{{\Bbb Z}}
\newcommand{\NN}{{\Bbb N}}
\newcommand{\QQ}{{\Bbb Q}}
\newcommand{\TT}{{\Bbb T}}
\newcommand{\GG}{{\Bbb G}}
\newcommand{\DD}{{\Bbb D}}
\newcommand{\NNN}{\NN\times \NN}
\newcommand{\demo}{\emph{Proof.} }
\newcommand{\Cd}{C_\bullet\,}
\newcommand{\Span}{\mathrm{span}\,}
\newcommand{\supp}{\mathrm{supp}\,}
\newcommand{\pa}{\prt_\alpha}
\newcommand{\pba}{\prt_{\sba}}
\newcommand{\be}{\begin{equation}}
\newcommand{\ee}{\end{equation}}
\newcommand{\nvar}[1]{\left\|#1\right\|_{var}}
\newcommand{\bmu}{{\mu_{inv}}}
\newcommand{\tov}{\mathop{\longrightarrow}\limits^{var}}
\newcommand{\eqd}{\mathop{=}\limits^{df}}
\newcommand{\indicator}[1]{1\!\!\mathrm I\{#1\}}
\newcommand{\brakes}[1]{\left(#1\right)}
\newcommand{\eqdistr}{\stackrel{\mathrm{d}}{=}}
\newcommand{\eqdef}{\stackrel{\mathrm{def}}{=}}
\newcommand{\limw}{\stackrel{\mathrm{w}}{\rightarrow}}
\newcommand{\lims}{\underset{n\rightarrow\infty}{{\lim\sup}}\,}
\newcommand{\norm}[1]{\left\Vert#1\right\Vert}
\newcommand{\abs}[1]{\left\vert#1\right\vert}
\newcommand{\figbrakes}[1]{\left\{#1\right\}}
\newcommand{\sqrbrakes}[1]{\left[#1\right]}
\newcommand{\mes}[1]{mes\{#1\}}
\newcommand{\Bc}{\mathcal{B}}
\newcommand{\Gc}{\mathcal{G}}
\newcommand{\Fc}{\mathcal{F}}
\newcommand{\Nc}{\mathcal{N}}
\newcommand{\Lc}{\mathcal{L}}


\title[Difference approximation for local times of multidimensional diffusions]
    {Difference approximation for local times of multidimensional diffusions}

\author{Alexey M. Kulik}
\address{Institute of Mathematics,
Ukrai\-ni\-an National Academy of Sciences, 3, Tereshchenkivska
Str., Kyiv 01601, Ukraine}
 \email{kulik@imath.kiev.ua}

\begin{abstract}
We consider sequences of additive functionals of difference
approximations for uniformly non-degenerate multidimensional
diffusions. The  conditions are given, sufficient for such a
sequence to converge weakly to a $W$-functional of the limiting
process. The class of the $W$-functionals, that can be obtained as
the limiting ones, is completely described in the terms of the
associated $W$-measures $\mu$  by the condition
 $$
\lim_{\delta\downarrow 0}\sup_{x\in\Re^m}\int_{\|y-x\|\leq
\delta}w(\|y-x\|)\mu(dy)=0,\quad\quad
w(r)\equiv \begin{cases} \max(-\ln r, 1),& m=2\\
r^{2-m},& m>2.
\end{cases}
$$
\end{abstract}



\thanks{The research was partially supported by the Ministry
of Education and Science of Ukraine, project N GP/F13/0095}

 \maketitle

\section{Introduction}

In the paper, we consider an $\Re^m$-valued diffusion process $X$
defined by an SDE
\be\label{00}X(t)=X(0)+\int_0^ta(X(s))\,ds+\int_0^tb(X(s))\,
dW(s), \quad t\in \ax, \ee and a sequence of processes $X_n, n\geq
1$, with their values  at the time moments ${k\over n}, k\in \NN$
given
 by a difference relation \be\label{01}X_n\brakes{k\over
n}=X_n\brakes{k-1\over n}+a\brakes{X_n\brakes{k-1\over n}}\cdot
{1\over n} +b\brakes{X_n\brakes{k-1\over n}}\cdot {\xi_k\over\sqrt
n}, \ee and, at all the other time moments, defined in a
piece-wise linear way:
 \be\label{02} X_n(t)=X_n\brakes{k-1\over
n}+(nt-k+1)\left[X_n\brakes{k\over n}-X_n\brakes{k-1\over n}
\right],\quad t\in\left[{k-1\over n},{k\over n}\right). \ee Here
and below,  $W$ is a Wiener process valued in $\Re^m$, $\{\xi_k\}$
is a sequence of i.i.d. random vectors in $\Re^m$, that belong to
the domain of attraction of the normal law, are centered and have
the identity for covariance matrix. Under standard assumptions
about coefficients of the equations (\ref{00}), (\ref{01}) (local
Lipschitz condition and linear growth condition), the
distributions of the processes $X_n$ in  $C(\ax, \Re^m)$ with the
given initial value $X_n(0)=x$ converge weakly to the distribution
of the process $X$ with $X(0)=x$ (see \cite{skor_as_m}).

In the paper, we deal with the following problem. Let
$\{\phi^{s,t}, 0\leq s\leq t\}$ be a $W$-functional of the
diffusion process defined by (\ref{00}). This, by definition,
means that (see \cite{dynkin}, Chapter 6) $\phi$ is a non-negative
 homogeneous additive functional with its characteristics
$$
\{f^{t}(x)\equiv E[\phi^{0,t}|X(0)=x], t\geq 0, x\in \Re^m\}$$
satisfying the condition $\sup_{x}f^{t}(x)<+\infty, t\in \ax$. We
consider a sequence of non-negative additive functionals
$\{\phi^{s,t}_n, 0\leq s\leq t\}, n\geq 1$ of the processes $X_n$,
of the form \be\label{03}\phi_{n}^{s,t}=\phi_{n}^{s,t}(X_n)\eqdef
{1\over n}\sum_{k:s\leq {k\over
n}<t}F_{n}\brakes{X_n\brakes{\frac{k}{n}}},\quad 0\leq s<t,\ee and
give the conditions sufficient for the joint distributions of
$(\phi_n,X_n),$ conditioned by $X_n(0)=x$ ($x\in\Re^m$ is an
arbitrary point), to converge weakly to the joint distribution of
$(\phi,X)$, conditioned by $X(0)=x$.

One motivation for posing such a problem is the following one. It
is well known that the theory of additive functionals of
$\Re^m$-valued Markov processes is closely related with the
potential theory. There exists a one-to-one correspondence between
$W$-functionals and so called $W$-measures (see \cite{dynkin},
Chapter 8); every  $W$-functional $\phi$ can be written at the
form  \be\label{04} \phi^{s,t}=\int_s^t{d\mu\over
d\lambda^m}(X(r))\,dr, \quad 0\leq s\leq t,  \ee where  $\mu$ is
the corresponding  $W$-measure, $\lambda^m$ is the Lebesgue
measure on $\Re^m$. In general, $W$-measure  $\mu$ is not
absolutely continuous w.r.t. $\lambda^m$; for singular $\mu$
equality (\ref{04}) is a formal notation, that can be
substantiated via an approximative procedure with $\mu$
approximated by an absolutely continuous measures. The functional
$\phi$, given by (\ref{04}), is called the local time of the
process $X$, corresponding  to the measure  $\mu$.

Given $W$-functional $\phi$, one can construct the sub-process
$X^\phi$ with its transition probability given by  \be\label{041}
P[X^\phi(t)\in \Gamma|X^\phi(0)=x]=E[\1_\Gamma(X(t))\cdot
\exp(-\phi^{0,t})|X(0)=x], \quad t\in \ax, x\in\Re^m, \Gamma\in
\Bf(\Re^m) \ee (see \cite{dynkin}, \S 6 of Introduction); the
measure $\mu$ is interpreted as the \emph{killing measure} for
this process. Let $X$ be the diffusion process given by the
equation (\ref{00}), then its generator is equal
$$
\Af=\sum_{i=1}^m a_i\cdot{\prt\over \prt x_i}+{1\over
2}\sum_{i,j=1}^m \sigma_{ij}\cdot {\prt^2\over \prt x_i\prt
x_j},\quad \sigma\equiv(\sigma_{ij})_{i,j=1}^m=bb^*.
$$
The well known Feynman-Kac formula gives the generator of the
process $X^\phi$ at the form  $\Af^\phi f=\Af f-{d\mu\over
d\lambda^m}\cdot f$, i.e., for every continuous bounded function
$g(\cdot),$  the solution to the following Cauchy problem for the
second order parabolic PDE \be\label{05} u'_t(t,x)=\sum_{i=1}^m
a_i\cdot{\prt\over \prt x_i}u (t,x)+{1\over 2}\sum_{i,j=1}^m
\sigma_{ij}\cdot {\prt^2\over \prt x_i\prt x_j}u (t,x)-{d\mu\over
d\lambda^m}(x) u(t,x), \quad u(0,x)=g(x) \ee has the following
probabilistic representation: \be\label{06} u(t,x)=E[g(X(t))\cdot
\exp(-\phi^{0,t})|X(0)=x].\ee Let us note that the term
${d\mu\over d\lambda^m}$ in the equation (\ref{05}), that
corresponds to the "heat flow-out", is a generalized function and
equation (\ref{05}) should be interpreted in the generalized sense
(therefore, it is natural to call this equation a singular one).
On the other hand, its solution is a bounded measurable function
due to the representation (\ref{06}).

Put \be\label{07} u_n(t,x)=E[g(X_n(t))\cdot
\exp(-\phi_n^{0,t})|X_n(0)=x].\ee The main result of the present
paper (Theorem  \ref{t21}) provides convergence in distribution of
the difference approximation $(\phi_n,X_n,)$ to $(\phi, X)$; thus,
under conditions of this theorem,
$$
u_n(t,x)\to u(t,x), \quad n\to \infty,\quad t\in \ax, x\in \Re^m.
$$
This means that this theorem, in particular, gives an opportunity
to apply, in a standard way, the Monte-Carlo method for numerical
solution of the Cauchy problem for the singular parabolic equation
(\ref{05}).

Let us make a short overview of the bibliography devoted to weak
convergence of the functionals of the type (\ref{03}). For
one-dimensional random walks and difference approximations of
one-dimensional diffusions, there exists a large variety of limit
theorems for the associated additive functionals. We do not
discuss these results in details, since we are mostly interested
in a multidimensional case, and refer to the monographs
\cite{skor_slob},\cite{bor_ibrag} and papers  \cite{gikh_1} --
\cite{Shi_Che_yor}.

In the multidimensional case the situation is essentially
different. The author does not know any paper where a limit
theorem for the functionals of the type (\ref{03}) would be proved
in the case where $X_n$ approximate a non-additive diffusion
process º $X$ (i.e, where the coefficients $a,b$ are
non-constant). The only multidimensional limit theorem, known for
the author, for the additive functionals of the type (\ref{03}),
is given in the paper \cite{bass_khoshn92} in the situation where
$X_n$ is a multidimensional aperiodic lattice random walk and $X$
is the Brownian motion in $\Re^m$ (also, the paper \cite{dynkin2}
deals with the closely related problems). The significant
difference between the results available in the one- and
multidimensional cases can be naturally explained by the fact that
the structure of the class of $W$-measures is much more
complicated in the second case than in the first one. For the
Brownian motion (and, also, for any non-degenerate diffusion), for
$m=1$, every  finite measure is a $W$-measure. For $m>1$, any
measure $\delta_z, z\in\Re^m$ is not a $W$-measure
($\Leftrightarrow$ there does not exists the local time at any
fixed point $\Leftrightarrow$ every one-point set has its capacity
equal to zero $\Leftrightarrow$ every one-point set is a polar
set). Therefore, in the case $m\geq 2$, the claims both on the
 "symbols"
${d\mu_n\over d\lambda^m}\equiv F_n$ of the approximating
aggregates (\ref{03}) and on the "symbol" ${d\mu\over d\lambda^m}$
of the limiting functional $\phi$ should be more delicate. In the
paper \cite{bass_khoshn92}, the uniform (w.r.t. $n$) analogue of
the following "dimensional" condition on the symbol ${d\mu\over
d\lambda^m}$ was used: \be\label{08} \exists\, C\in\ax,
\gamma>m-2:\quad \mu(\{y|\|y-x\|\leq r\})\leq Cr^\gamma,\quad
r\geq 0. \ee For the Brownian motion in $\Re^m$, the following
criterium is well known (\cite{dynkin}, Chapter 8): measure $\mu$
is a $W$-measure iff \be\label{09}
\sup_{x\in\Re^m}\int_{\|y-x\|\leq 1}w(\|y-x\|)\mu(dy)<+\infty,\ee
$$
w(r)\equiv \begin{cases} \max(-\ln r, 1),& m=2\\
r^{2-m},& m>2.
\end{cases}
$$
It is easy to verify that the condition (\ref{08}) is sufficient
for the condition (\ref{09}) to hold true. However, it is not a
necessary one (see Example \ref{e51}). This, in particular, means
that in the main  limit theorem of \cite{bass_khoshn92} only the
functionals from some proper subclass of the class of
$W$-functionals (namely, the functionals with their $W$-measures
satisfying the "dimensional" condition (\ref{08})) can be obtained
as a limiting ones.

The main result of the present paper (Theorem \ref{t21})
establishes the weak convergence of the functionals (\ref{03}) for
difference approximations $X_n$ of multidimensional uniformly
non-degenerate diffusions $X$. In our framework, the class of the
difference approximations $X_n$ is wide enough. We claim the
densities of the transition probabilities for  $X_n$ to satisfy a
proper version of the local limit theorem  (property  B4), Chapter
4 below). We rely on the results of the papers \cite{kon_mammen},
\cite{kon} while giving conditions, sufficient for such a claim.

The condition, imposed in Theorem \ref{t21} on the "symbols" of
the approximating aggregates, is the uniform (w.r.t. $n$) analogue
of the condition \be\label{010} \lim_{\delta\downarrow
0}\sup_{x\in\Re^m}\int_{\|y-x\|\leq \delta}w(\|y-x\|)\mu(dy)=0.
\ee Condition (\ref{08}) is sufficient, but not necessary one for
the condition (\ref{010}) to hold true (see Example \ref{e51}).
Let us discuss the relation between conditions (\ref{010}) and
(\ref{09}) in a more details. In the present paper, in order to
make exposition more short and transparent, we consider the
measures $\mu$ with a compact supports, only. For such a measure,
using the standard estimate for the transition density of a
non-degenerate diffusion (see \cite{dynkin}, Appendix, \S 6 and
references there), one can check, in a standard way, that the
condition (\ref{09}) is also the necessary and sufficient
condition for a measure $\mu$ to be a $W$-measure for uniformly
non-degenerate diffusion. The condition (\ref{010}) is clarified
by the following statement.

\begin{prop}\label{p01} Let $X$ be a  diffusion process, valued in $\Re^m, m\geq
2$,with its coefficients satisfying   condition A1) of Theorem
\ref{t21}, given below. Let  $\mu$ be a $W$-measure with a compact
support, $\phi$ be the $W$-functional, corresponding to this
measure and $f$ be its characteristics. Then the following
statements are equivalent:

(i) for any  $t\in\ax,$ the function $f^t(\cdot)$ is uniformly
continuous on $\Re^m$;

(ii) $\sup_{x\in\Re^m}f^\delta(x)\to 0,\, \delta\downarrow 0$;

(ii³) $\mu$ satisfies the condition (\ref{010});
\end{prop}

We prove Proposition \ref{p01} in the Chapter 4. One can interpret
the statement of this Proposition in the following way: any
$W$-measure, satisfying (\ref{010}), correspond to the
$W$-functional that is regular w.r.t. the phase variable. The
class of functionals, that can be obtained as the limit ones in
the context of the main result the present paper (Theorem
\ref{t21}), exactly coincides with the class of the functionals,
regular w.r.t. the phase variable in a sense given by Proposition
\ref{p01} (see Remark \ref{r23}).

\section{The main statement}

We consider the objects, given by (\ref{00}) -- (\ref{03}), for
$m\geq 2$. We use notation $\|\cdot\|$ for the Euclidean norm, not
indicating explicitly the space this norm is written for. The
classes of functions, that have  $k$ continuous derivatives, and
functions, that are continuous and bounded together with their $k$
derivatives, are denoted by $C^k$ and $C^k_b$, correspondingly.
The derivative (the gradient) is denoted by $\nabla$. The weak
convergence of the (not necessary probability) measures $\mu_n$ to
$\mu$ means, by the definition, convergence
$\int_{\Re^m}f\,d\mu_n\to \int_{\Re^m}f \, d\mu_n$ for every $f\in
C_b(\Re^m)$.

Let us formulate the main conditions on the objects involved to
(\ref{00}) -- (\ref{03}).

\begin{itemize} \item[A1)] $a\in C_b^2(\Re^m,\Re^m), b\in C^2_b(\Re^m,\Re^{m\times m})$
and there exist positive constants  $c,C$ such that
$$
c\|\theta\|^2\leq (b(x)b^*(x)\theta, \theta)_{\Re^m}\leq
C\|\theta\|^2, \quad x,\theta\in \Re^m.
$$
Furthermore, the function $\nabla^2 b$ satisfies the H\"older
condition with some positive exponent.

\item[A2)] I.i.d. random vectors $\{\xi_k\}$ are centered and
have the identity for covariance matrix.

\item[A3)] Random vectors $\{\xi_k\}$ possess the distribution
density $p\in C^4(\Re^m)$. There exists a function $\psi:\Re^m\to
\ax$ such that
$$
\sup_{x\in\Re^m}\psi(x)<+\infty, \quad
\int_{\Re^m}\|x\|^{m^2+2m+4}\psi(x)\, dx<+\infty
$$
and
$$
\|[\nabla^i p](x)\|\leq \psi(x), \quad x\in \Re^m, \quad
i=0,\dots, 4.
$$
\item[A4)] $F_n(x)\geq 0, x\in\Re^m, n\geq 1$ and ${1\over
n}\sup\limits_{x\in\Re^m}F_n(x)\to 0, n\to\infty$.

\item[A5)] Measures $\mu_n(dx)\equiv F_n(x)\lambda^m(dx)$ weakly
converge to the finite measure $\mu$, that has a compact support.

\item[A6)]  The following uniform analogue of the condition
(\ref{010}) holds true:
$$
\lim_{\delta\downarrow
0}\mathop{\lim\sup}_{n\to+\infty}\sup_{x\in\Re^m}\int_{\|y-x\|\leq
\delta}w(\|y-x\|)\mu_n(dy)\to 0.
$$
\end{itemize}

\begin{zau} If  conditions A5) and
A6) hold true, then the measure $\mu$  satisfies condition
(\ref{010}). In particular, $\mu$ is a $W$-measure.
\end{zau}
\begin{zau}\label{r22} The function $w(\cdot)$ is bounded on $[\delta,+\infty)$
for any given $\delta>0$. Thus, conditions A5) and A6) imply that
$\mathop{\lim\sup}_n \mu_n(\Re^m)<+\infty$ and
$$
\mathop{\lim\sup}_{n\to+\infty}\sup_{x\in\Re^m}\int_{\Re^m}w(\|y-x\|)\mu_n(dy)<+\infty.
$$
\end{zau}

Let us proceed with the formulation of the main statement.
Together with the functionals $\phi_n$, that are piece-wise
constant w.r.t. both time variables, we consider the random broken
lines, constructed from these functionals:
$$
\psi_n^{s,t}=\phi_{n}^{{j-1\over n}, {k-1\over
n}}+(ns-j+1)\phi_{n}^{{j-1\over n}, {j\over
n}}+(nt-k+1)\phi_{n}^{{k-1\over n}, {k\over n}},\,\,
s\in\left[{j-1\over n},{j\over n}\right),t\in\left[{k-1\over
n},{k\over n}\right).
$$
We interpret the random broken lines  $\psi_n$ as the random
elements, taking values in $C(\TT,\Re^+)$, where
$\TT\eqdef\{(s,t)|0\leq s\leq t\}$.  The $W$-functional
$\phi=\phi(X)$ of the process $X$ we define by the formula
(\ref{04}) (the measure $\mu$ is taken from the condition A5)).

\begin{thm} \label{t21} Let the conditions  A1) -- A6) hold.
Then $(X_n,\psi_n(X_n))\Rightarrow (X,\phi(X))$ in a sense of weak
convergence in $C(\ax, \Re^m)\times C(\TT,\Re^+)$.
\end{thm}

We prove Theorem \ref{t21} in the Chapter 4.

\begin{zau}\label{r23} Let an arbitrary  $W$-measure
$\mu$, satisfying condition (\ref{010}), be given. Then one can
construct a sequence of the functions $\{F_n\}$ in such a way that
conditions A4)-A6) hold true. For instance, one can define $F_n$
by
$$
F_n(x)=n\cdot\int_0^{1\over n}\int_{\Re^m}q_{t}(x,y)\mu(dy)\,
dt,\quad x\in\Re^m,
$$
where  $q_t(x,y)=(2\pi t)^{-{m\over
2}}\exp\brakes{-{\|y-x\|^2\over 2}}$ is the transition probability
density for the Brownian motion  $W$. Then ${1\over n}F_n$ is
equal to the value at the point ${1\over n}$ of the
characteristics of the $W$-functional of the process  $W$,
$$
L_t=\int_0^t{d\mu\over d\lambda^m}(W_s)\,ds,
$$
and the properties A4)-A6) can be proved analogously to the proof
of Proposition \ref{p01} (see Chapter 4).
\end{zau}

\section{Weak convergence of additive functionals of
a sequence of Markov chains}

Our proof of Theorem  \ref{t21} is based on the general theorem on
convergence in distribution of a sequence of additive functionals
of Markov chains, given in the paper \cite{kar_kul}. In this
chapter, we give a detailed exposition of the objects and
auxiliary notions, that are used in this theorem.

In this chapter we suppose that the processes  $X_n(\cdot),
X(\cdot)$ are defined $\ax$ and take their values in a locally
compact metric space $(\Xf,\rho)$. We say that the process $X$
possesses the Markov property at the time moment $s\in\ax$ w.r.t.
filtration  $\{\mathcal{G}_t,t\in\ax\}$, if $X$ is adapted with
this filtration and for every  $k\in \NN, t_1,\dots,t_k>s$ there
exists a probability kernel  $\{P_{st_1\dots t_k}(x,A), x\in \Xf,
A\in\mathcal{B}(\Xf^k)\}$ such that  \be\label{30}
E[\1_{A}((X(t_1),\dots,X(t_k)))|\mathcal{G}_s]=P_{st_1\dots
t_k}(X(s),A) \quad \hbox{a.s.,} \quad A\in\mathcal{B}(\Xf^k). \ee
The measure $P_{st_1\dots t_k}(x,\cdot)$ has a natural
interpretation as the conditional finite-dimensional distribution
of $X$ at the points $t_1,\dots,t_k$ under condition $\{X(s)=x\}$;
below, we use notation $P_{st_1\dots
t_k}(x,\cdot)=P((X(t_1),\dots,X(t_k))\in\cdot|X(s)=x)$.

Everywhere below we claim the process $X$ to possess the Markov
property w.r.t. its canonical filtration at every $s\in \ax$, and
every processes $X_n$ to possess this property (w.r.t. its
canonical filtrations) at the points of the type ${i\over n}, i\in
\ZZ_{+}$; this means that every process $X_n$ is, in fact, a
Markov chain with the time scale, proportional to  ${1\over n}$.

Let the additive functionals $\phi_n$ be given by the formula
(\ref{03}). For the functional $\phi_n$, its characteristics $f_n$
(the analogue of the characteristics of a $W$-functional) is
defined by the formula \be\label{32} f_n^{s,t}(x) \eqdef
E[\phi_{n}^{s,t}(X_n)|X_n(s)=x],\quad s={i\over n}, i\in \ZZ_+,
t>s, x\in \Xf. \ee Let us note that the process $X_n$ possesses
the Markov property w.r.t. its canonical filtration at the points
$s={i\over n}, i\in \ZZ_+$ and the functional (\ref{03}) is the
function of the values of $X_n$ at the finite family of teh time
moments. Therefore the mean value in (\ref{32}) is well defined as
the integral w.r.t. family of the conditional finite-dimensional
distributions  $\{P_{st_1\dots t_k}(x,\cdot), t_1,\dots,t_k>s,
k\in \NN\}$ of the process $X_n$.

The following result (\cite{kar_kul}, Theorem 1) is an analogue of
the well known theorem by E.B.Dynkin, that describes convergence
of  $W$-functionals in the terms of their characteristics
 (\cite{dynkin}, Theorem 6.4).

\begin{thm}\label{t31} Let the sequence of the processes $X_n$ be given, providing Markov
approximation for the homogeneous Markov process  $X$ (see
Definition \ref{d31} below), and let the sequence
  $\{\phi_n\equiv \phi_n(X_n)\}$ be defined by
(\ref{03}).  Suppose that the following conditions hold true:
\begin{enumerate}
\item The functions ${1\over n}F_{n}(\cdot)$ are non-negative,
bounded on $\Xf$ and uniformly converge to zero:
$$
\delta(F_{n})\eqdef {1\over n}\sup_{x\in\Xf}F_{n}(x)\to 0,\quad
n\to \infty.
$$
 \item There exists a function $f$, that is a characteristics of a certain
 $W$-functional  $\phi=\phi(X)$ of the limiting process $X$, such that, for every
  $T\in\ax$,
$$\underset{s={i\over n}, t\in (s,T)}{\sup}\sup_{x\in \Xf}|
f^{s,t}_n(x)-f^{t-s}(x)|\rightarrow0, \quad n\rightarrow\infty.$$

\item  The limiting function $f$ is continuous w.r.t. variable
$x$, locally uniformly w.r.t. time variable, i.e., for every
  $T\in\ax$,
$$\sup\limits_{t\leq T}\abs{f^{t}(x)-f^{t}(y)}\rightarrow0,\quad \|x-y\|\rightarrow0.$$
\end{enumerate}
Then, for the random broken lines $\psi_n$, corresponding to
$\phi_n$,
$$\psi_n(X_n)\Rightarrow\phi(X)\equiv \{\phi^{s,t}(X), (s,t)\in
\TT\}$$ in a sense of weak convergence in $C(\TT,\Re^+)$.

If, additionally,  $X_n\Rightarrow X$ ó  $C(\ax, \Xf)$, then
 $(X_n,\psi_n(X_n))\Rightarrow
(X,\phi(X))$ n a sense of weak convergence in $C(\ax, \Xf)\times
C(\TT,\Re^+)$.
\end{thm}

The notion of Markov approximation, introduced in  \cite{kulik},
is a key one in Theorem  \ref{t31}.

\begin{ozn}\label{d31}  The sequence of the processes $\{X_n\}$ provides
 the Markov approximation for the
Markov process $X$, if for every $\gamma>0, T<+\infty$ there exist
a constant $K(\gamma,T)\in \NN$ and a sequence of two-component
processes $\{\hat Y_n=(\hat X_n,\hat X^n)\}$, possibly defined on
another probability space, such that

(i) $\hat X_n\mathop{=}\limits^d X_n, \hat X^n \mathop{=}\limits^d
X$;

(ii) the processes $\hat Y_n$,$\hat X_n, \hat X^n$ possess the
Markov property at the points ${i K(\gamma,T)\over n}, i\in \NN$
w.r.t. the filtration $\{\hat {\Fc}_t^n=\sigma(\hat Y_n(s),s\leq
t)\};$

(iii)  $\lim\sup_{n\to+\infty}P\left(\sup_{i\leq {Tn\over
K(\gamma,T)}} \rho\left(\hat X_n\left({i K(\gamma,T)\over
n}\right),\hat X^n\left({i K(\gamma,T)\over
n}\right)\right)>\gamma\right)<\gamma.$
\end{ozn}

The following result, on the one hand, provides an example, that
clarifies the given above definition, and, on the other hand,
gives the opportunity to apply Theorem \ref{t31} in order to prove
the main statement of the paper.

\begin{lem}\label{l31} Let $a,b$ be bounded and satisfy Lipschitz condition,
 $m\geq 1$ be arbitrary and i.i.d. random vectors $\{\xi_k\}$ be centered,
have  the identity for covariance matrix and satisfy condition
$E\|\xi_k\|^{2+\delta}<+\infty$ for some $\delta>0$.

Then the sequence of the processes $X_n$, given by
(\ref{01}),(\ref{02}), provides the Markov approximation for the
process $X$, given by the equation (\ref{00}).
\end{lem}

At the Example 3 of the paper \cite{kar_kul}, the statement of the
Lemma was proved with the use of the pathwise uniqueness property
of the equation (\ref{00}). The reasonings of such a kind are a
qualitative ones, and can not provide explicit estimates for the
rate of convergence. Therefore, here we give another
straightforward proof, that gives possibility for the further
estimates and generalizations.

\demo We start from the construction, described in the proof of
Theorem 1 \cite{kulik}. Denote $S_N=\sum_{k=1}^N \xi_k$; due to
CLT, $n^{-{1\over 2}}S_n\Rightarrow W(1)$.  Condition
$E\|\xi_k\|^{2+\delta}<+\infty$ ensures that the family
$\{{S_N^2\over N}\}$ is uniformly integrable, and therefore the
Wasserstein distance between the distributions of the vectors
$N^{-{1\over 2}}S_N$ and $W(1)$ tends to 0 as $N\to \infty$. This
means that, for any $\eps>0$, there exist  $N_\eps\in \NN$ and
random vector $(\eta_\eps, \zeta_\eps)$ such that
$$
E\|\eta_\eps-\zeta_\eps\|^2_{\Re^m}<\eps, \quad \eta_\eps\eqd
{S_{N_\eps}\over \sqrt{N_\eps}},\quad \zeta_\eps\eqd W(1).
$$
Let $\eps>0$ be fixed; we construct the probability space
$(\Omega^1,\Fc^1,P^1)$ in the following way:
$\Omega^1=(\Re^m)^{N_\eps}\times C([0,1])$, $\Fc^1=\Bc(\Omega^1)$.
Denote the coordinates of a point  $\omega^1\in\Omega^1$ by
$\chi=(\chi_1,\dots,\chi_{n_\eps})\in (\Re^m)^{N_\eps}, \phi\in
C([0,1])$. Define the following measures: $Q(du,dv)$ is the joint
distribution of  $(\eta_\eps,\zeta_\eps)$, $U_\eps(d\chi,u)$ is
the conditional distribution of  $\{\xi_1,\dots,\xi_{N_\eps}\}$
under condition $\{{S_{N_\eps}\over \sqrt{N_\eps}}=u\}$, and
$V_\eps(d\phi,v)$ is the conditional distribution of $W(\cdot)$
under condition $\{W(1)=v\}$. We put
$$
P^1(A)=\int_{(\Re^m)^2}\left[\int_{A}
U_\eps(d\chi,u)V_\eps(d\phi,v)\right]Q(du,dv),\quad A\in {\Fc}^1.
$$
Now we define the probability space  $(\Omega,\Fc,P)$ as the
infinite product of the copies of $(\Omega^1,\Fc^1,P^1)$. For
$\omega=(\chi^1,\phi^1,\chi^2,\phi^2,\dots)\in \Omega$ define the
sequence  $\{\hat \xi_k(\omega), k\geq 1\}$ by the formula
$$
\hat \xi_1(\omega)=\chi^1_1,\,\hat
\xi_2(\omega)=\chi_2^1,\,\dots,\, \hat
\xi_{N_\eps}(\omega)=\chi^1_{n_\eps},\,\xi_{N_\eps+1}(\omega)=\chi^2_1,\,\dots,
$$
and the process  $\{\hat W^n(t),t\in\ax\}$ by the formula
$$
\hat W^n(t)(\omega)={1\over \sqrt
{n}}\left[\sum_{k=1}^{[nt]}\phi^k(1)+\phi^{[nt]+1}\Bigl({t-[nt]\over
n}\Bigr)\right], \quad t\geq 0.
$$
By the construction, the sequence $\{\hat \xi_k\}$ has the same
distribution with the sequence $\{\xi_k\}$ and the process $\hat
W^n$ is a Brownian motion in $\Re^m$. Now, let us define processes
$\hat X_n, \hat X^n$ by the formulae (\ref{01}),(\ref{02}) and
(\ref{00}), with $\{\xi_k\}$ replaced by  $\{\hat \xi_k\}$ and $W$
replaced by $\hat W^n$; by the construction, the process $\hat
Y_n=(\hat X_n,\hat X^n)$ satisfies the condition (i) of Definition
\ref{d31}.

Also, by the construction, the sets
$$
\Xi_n^l=\left\{\hat \xi_{l N_\eps+1}, \dots, \hat
\xi_{(l+1)N_\eps}, \brakes{\hat W^n(\cdot)-\hat W^n\brakes{l
N_\eps\over n}}\Big|_{[{l N_\eps\over n}, {(l+1)N_\eps\over
n}]}\right\}, \quad l=0,1,\dots
$$
are mutually independent. According to (\ref{01}),(\ref{02}) and
(\ref{00}), the value of the process $\hat Y_n$ at the given time
moment ${i N_\eps\over n}, i\in \NN,$ is a functional of
$\Xi_n^0,\dots, \Xi_n^{i-1}$, and the values of  $\hat Y_n, \hat
X_n$ or $\hat X^n$, at any time moment $t>{i N_\eps\over n}$, are
a functionals of $\Xi_n^{i},\Xi_n^{i+1},\dots$ and $\hat
Y_n\brakes{i N_\eps\over n}, \hat X_n\brakes{i N_\eps\over n} $ or
$\hat X^n\brakes{i N_\eps\over n}$, respectively. Thus, the
processes $\hat Y_n, \hat X_n$ and $\hat X^n$ possess the Markov
property w.r.t. filtration $\{\hat \Fc_t, t\in \ax\}$, generated
by $\hat Y_n$, at the time moments ${i N_\eps\over n}, i\in \NN$.

Let us proceed with the estimation of the distance between  $\hat
X_n$ and $\hat X^n$. In order to shorten exposition, we will give
the estimate in the partial case $m=1, a\equiv 0$; in general
case, the argumentation is completely analogous, but the
calculations take  more place. Denote
$$
t_i={i N_\eps\over n},  \quad \Delta_{i,n}=\hat X_n(t_i)-\hat
X_n(t_{i-1})=\sum_{l=1}^{N_\eps} b\brakes{\hat
X_n\brakes{t_{i-1}+{l-1\over n}}}\cdot {\hat \xi_{(i-1)
N_\eps+l}\over \sqrt{n}},
 $$
$$
 \Delta_i^n=\hat X^n(t_i)-\hat X^n(t_{i-1})=\int_{t_{i-1}}^{t_i}
b\brakes{\hat X^n(s)}\, dW^n(s),\quad   i\in \NN.
$$
By the construction, $E[\Delta_{i,n}|\hat
\Fc_{t_{i-1}}]=E[\Delta_i^n|\hat \Fc_{t_{i-1}}]=0$, therefore,
\be\label{321} E(\hat X_n(t_i)-\hat X^n(t_i))^2=E(\hat
X_n(t_{i-1})-\hat X^n(t_{i-1}))^2+E(\Delta_{i,n}-\Delta_i^n)^2,
\quad i\in \NN. \ee We write the  decomposition
$$
\Delta_{i,n}=b\brakes{\hat X_n\brakes{t_{i-1}}}\cdot
\sum_{l=1}^{N_\eps}  {\hat \xi_{(i-1) N_\eps+l}\over
\sqrt{n}}+\Gamma_{i,n},
$$
$$
\Gamma_{i,n}\equiv \sum_{l=1}^{N_\eps} \left[b\brakes{\hat
X_n\brakes{t_{i-1}+{l-1\over n}}}-b\brakes{\hat
X_n\brakes{t_{i-1}}}\right]\cdot {\hat \xi_{(i-1) N_\eps+l}\over
\sqrt{n}}.
$$ The pair of the processes $(\{\hat X_n(t), t\in \ax\}, \{\hat \xi_k, k\in
\NN\})$ has the same distribution with the pair  $(\{X_n(t), t\in
\ax\}, \{\xi_k, k\in \NN\})$. Thus, for every $k\geq 1$, the
random variable  $\hat \xi_k$ does not depend on the values of the
process $\hat X_n$ on the interval $[0, {k-1\over n}]$. Using
this, and taking into account that the function $b$ is bounded by
a constant $B$ and satisfies the Lipschitz condition with a
constant $L$, we get the estimate  \be\label{33} E
\Gamma_{i,n}^2\leq {L^2\over n}\sum_{l=1}^{N_\eps}E\left[\hat
X_n\brakes{t_{i-1}+{l-1\over n}}-\hat
X_n\brakes{t_{i-1}}\right]^2\leq {L^2\over
n}\sum_{l=1}^{N_\eps}B^2\cdot {l-1\over n}<{L^2B^2N_\eps^2\over
2n^2}. \ee Let us also write  the decomposition
$$
\Delta_i^n=b\brakes{\hat X^n\brakes{t_{i-1}}}\cdot \Big[\hat
W^n(t_i)-\hat W^n(t_{i-1})\Big]+\Gamma_i^n,
$$
where the second summand (we do not write it explicitly) can be
estimated analogously to  (\ref{33}): \be\label{34} E
\left[\Gamma_i^n\right]^2\leq {L^2B^2N_\eps^2\over 2n^2}. \ee At
last, we write the decomposition
$$
b\brakes{\hat X_n\brakes{t_{i-1}}}\cdot \sum_{l=1}^{N_\eps}  {\hat
\xi_{(i-1) N_\eps+l}\over \sqrt{n}}-b\brakes{\hat
X^n\brakes{t_{i-1}}}\cdot \Big[\hat W^n(t_i)-\hat
W^n(t_{i-1})\Big]=
$$
$$
=\left[b\brakes{\hat X_n\brakes{t_{i-1}}}-b\brakes{\hat
X^n\brakes{t_{i-1}}}\right]\cdot \Big[\hat W^n(t_i)-\hat
W^n(t_{i-1})\Big]+\sqrt{N_\eps\over n}\cdot b\brakes{\hat
X_n\brakes{t_{i-1}}}\cdot \Upsilon_n^i,
$$
where  $\Upsilon_n^i={1\over \sqrt {N_\eps}}\sum_{l=1}^{N_\eps}
\hat \xi_{(i-1) N_\eps+l}-[t_i-t_{i-1}]^{-{1\over 2}}\Big[\hat
W^n(t_i)-\hat W^n(t_{i-1})\Big]$. By the construction, for every
$i\geq 1,$ the random variables  $\Big[\hat W^n(t_i)-\hat
W^n(t_{i-1})\Big]$ and  $\Upsilon_n^i$ do not depend on the values
of the process $\hat Y_n$ on the interval $[0, {(i-1)Ò_\eps\over
n}]$. Moreover, $E[\Upsilon_n^i]^2<\eps.$ Therefore, the following
estimates hold true \be\label{35} E\left[b\brakes{\hat
X_n\brakes{t_{i-1}}}-b\brakes{\hat X^n\brakes{t_{i-1}}}\right]^2
\Big[\hat W^n(t_i)-\hat W^n(t_{i-1})\Big]^2\leq {L^2N_\eps\over
n}E(\hat X_n(t_{i-1})-\hat X^n(t_{i-1}))^2,\ee
\be\label{36}E\left[b\brakes{\hat X_n\brakes{t_{i-1}}}\cdot
\Upsilon_n^i\right]^2\leq B^2\eps^2. \ee Now, using the
decomposition (\ref{321}), the Cauchy inequality  $(\sum_{k\leq 4}
x_k)^2\leq 4\sum_{k\leq 4} x_k^2$ and the estimates
(\ref{33})--(\ref{36}), we obtain  \be\label{37} E(\hat
X_n(t_i)-\hat X^n(t_i))^2\leq \left(1+{4L^2N_\eps \over
n}\right)E(\hat X_n(t_{i-1})-\hat
X^n(t_{i-1}))^2+{4L^2B^2N_\eps^2\over n^2}+{4B^2\eps N_\eps\over
n}. \ee Iterating (\ref{37}), we obtain \be\label{38} E(\hat
X_n(t_i)-\hat X^n(t_i))^2\leq \left[{4L^2B^2N_\eps\over
n}+4B^2\eps\right]\cdot {N_\eps\over n}\cdot\sum_{j\leq
i}\left(1+{4L^2 N_\eps\over n}\right)^j.\ee Let  $t_i\leq T$, that
means that $i\leq {Tn\over N_\eps}$. Then the sum in the right
hand side of (\ref{38}) contains at most  ${Tn\over N_\eps}$
summands, and every summand is not greater than $\exp\left[{4L^2
N_\eps\over n}\cdot{Tn\over N_\eps}\right]=e^{4L^2 T}. $ This
provides the estimate  \be\label{39} E(\hat X_n(t_i)-\hat
X^n(t_i))^2\leq \left[{4L^2B^2N_\eps T\over
n}+4B^2T\eps\right]\cdot e^{4L^2 T}.\ee The sequence  $\{(\hat
X_n(t_i)-\hat X^n(t_i)), i\in \NN\}$, by the construction, is a
martingale, thus, using the maximal martingale inequality
(\cite{Doob}, Ch. VII, \S 3), we obtain the estimate
\be\label{310} P\left[\max_{i:t_i\leq T}(\hat X_n(t_i)-\hat
X^n(t_i))^2\geq \gamma^2\right]\leq
\gamma^{-2}\left[{4L^2B^2N_\eps T\over n}+4B^2T\eps\right]\cdot
e^{4L^2T}, \quad \gamma>0.\ee Now we can complete the proof of the
Lemma. For a given $\gamma, T$ choose $\eps>0$ in such a way that
$$
16B^2T\cdot e^{4L^2T}\cdot \eps<\gamma^3,
$$
and proceed with the construction, described above, with this
$\eps$. We have already seen that, under this construction,
conditions (i) and (ii) hold true with $K(\gamma, T)=N_\eps$. The
estimate (\ref{310}) provides that the condition (iii) holds true
with the same $K(\gamma, T)$. The lemma is proved.

\begin{zau} Denote by $\mathbf{K}(\gamma, T)$ the minimum of the
set of such numbers $K\in \NN$, that here exists a process $\hat
Y_n,$ satisfying conditions (i) -- (iii) of Definition  \ref{d31}
with $K(\gamma, T)=K$. In the paper \cite{kulik}, it is shown (the
part II of Theorem 1) that, in the basic case  $a\equiv 0, b\equiv
I_{\Re^m}$,
$$
\sup_{\gamma>0}\mathbf K(\gamma, T)=+\infty
$$
as soon as the distribution of $\xi_1$ differs from the normal
one. One can say that, while the accuracy of the approximation
becomes better (the accuracy is described by the parameter
$\gamma$), the Markov properties  of the two-component process
necessarily become worse (these properties are described by
$\mathbf{K}(\gamma,T)$).
\end{zau}

\section{The proofs of Theorem \ref{t21} and Proposition \ref{p01}}

We reduce the proof of Theorem \ref{t21} to the verification of
the conditions of Theorem \ref{t31}. The sequence $X_n$ provides
the Markov approximation for the process $X$ due to Lemma
\ref{l31}. Condition 1 of Theorem \ref{t31} holds true due to
condition  A4). Let us check that the conditions 2 and 3 hold
true. The characteristics $f^t$ of the functional
$\phi^{s,t}=\int_s^t{d\mu\over d\lambda^m}\brakes{X(r)}\, dr$ has
the form $$ f^t(x)=\int_0^t\int_{\Re^m}p_r(x,y)\mu(dy)\,dr,\quad
x\in\Re^m, t\geq 0,$$  where  $\{p_t(x,y), t\geq 0, x,y\in
\Re^m\}$ is the transition probability density for the process
$X$. Existence of such a density under condition A1) is a standard
result of the theory of parabolic equations. Moreover, this
density possesses the following properties   (see, for instance,
 Appendix to \cite{dynkin},\S 6, and references there).

\begin{itemize}
\item[B1)] The function $p:(t,x,y)\mapsto p_t(x,y)$ is uniformly
continuous on $[\delta,T]\times\Re^m\times\Re^m$ for every
$0<\delta<T$.

\item[B2)] There exist a constants  $M,\alpha>0$ such that
$$
p_t(x,y)\leq Mt^{-{m\over 2}}\exp\brakes{-{\alpha\|y-x\|^2\over
t}}.
$$
\item[B3)] There exist a constants $M_1,
M_2,\alpha,\beta,\lambda>0$ such that
$$
p_t(x,y)\geq M_1t^{-{m\over 2}}\exp\brakes{-{\alpha\|y-x\|^2\over
t}}-M_2t^{-{m\over 2}+\lambda}\exp\brakes{-{\beta\|y-x\|^2\over
t}}.
$$
\end{itemize}

Under conditions A2),A3), the processes $X_n$ possess the
transition probability densities $p_t^n(x,y)$ at the time moments
$t\in {1\over n}\NN$, that means that
$$
P(X_n(t)\in \Gamma|X_n(s)=x)=\int_\Gamma p^n_{t-s}(x,y)\,dy, \quad
s,t\in {1\over n}\ZZ, s<t, \quad x\in \Re^m, \Gamma\in \Bf(\Re^m),
$$
and, moreover, for every $t\in {1\over n}\NN$, the function
$p_t^n$ is a continuous one. The characteristics $f_n$ of the
functionals $\phi_n$ can be expressed through these densities by
the formula  \be \label{40} f_{n}^{s,t}(x)={1\over
n}F_n(x)+{1\over n}\sum_{k\in \NN, {k\over
n}<t-s}\int_{\Re^m}F_{n}(y) p_{k\over n}^{n}(x,y)\,
dy=f_n^{0,t-s}(x), \quad s\leq t, x\in\Re^m. \ee Below, we denote
$f_n^{0, t}=f_n^t$.  Theorem 2.1 \cite{kon_mammen} and Theorem 1
\cite{kon} imply that, under conditions A1) -- A3), the following
estimate for the deviation of the densities $p_t^n$ from the
limiting density $p_t$ holds true.
\begin{itemize}
\item[B4)] For  any  $T\in \ax$,
$$\sup_{n\in\NN, t\leq T}\sup_{x,y\in \Re^m}\sqrt{n}\cdot t^{-{m\over
2}}\cdot \brakes{1+\brakes{\|y-x\|\over \sqrt{t}}^m}\cdot
|p_t(x,y)-p_t^n(x,y)|<+\infty.
$$
\end{itemize}

Denote, for  $\delta>0$,
$$
f^t_\delta(x)=f^t(x)-f^{\delta\wedge t}(x)=\int_{\delta\wedge
t}^t\int_{\Re^m}p_r(x,y)\mu(dy)\,dr,\quad x\in\Re^m, t\geq 0,
$$
$$
f_{n,\delta}^{t}(x)=f^t_n(x)-f^{\delta\wedge t}_n(x)={1\over
n}\sum_{k\in \NN, \delta\leq {k\over n}<t}\int_{\Re^m}F_{n}(y)
p_{k\over n}^{n}(x,y)\, dy, \quad s\leq t, x\in\Re^m.
$$
Conditions  A4), A5) and properties B1),B2),B4) imply the
following statement.

\begin{prop}\label{p41} For any $\delta>0$,

(i) the function $f_\delta:(t,x)\mapsto f^t_\delta(x)$ is
uniformly continuous on $[0,T]\times \Re^m$;

(ii) the functions  $f_{\delta,n}:(t,x)\mapsto f^t_{n,\delta}(x)$
converge, as $n\to \infty$, to the function $f_{\delta}$ uniformly
on  $[0,T]\times \Re^m$.
\end{prop}

The proof of Proposition \ref{p41} is quite standard and we omit
it here (see, for instance, the arguments given in the proof of
Theorem 3 \cite{kar_kul}). It follows from Proposition \ref{p41}
that, in order to prove that the conditions 2 and 3 of Theorem
\ref{t31} hold true (and therefore, to prove the Theorem
\ref{t21}), it is sufficient to prove the following two relations:
\be\label{41}\lim_{\delta \downarrow 0} \sup_{x\in\Re^m}
f^\delta(x)=0, \ee \be\label{42}\lim_{\delta \downarrow 0}
\mathop{\lim\sup}_{n\to+\infty}\sup_{x\in\Re^m} f_n^\delta(x)=0.
\ee Let us prove (\ref{42}). In the exposition below, we suppose
that  $\delta\leq 1$. It follows from B2),B4) that, for any
$T\in\ax$, there exists a constant $C_T\in \ax$ such that
$$
p_t^n(x,y)\leq C_T\cdot t^{-{m\over
2}}\left[\exp\brakes{-{\alpha\|y-x\|^2\over
t}}+\brakes{1+\brakes{\|y-x\|\over \sqrt{t}}^m}^{-1}\right],\quad
t\in {1\over n}\NN,t\leq T,  n\in \NN.
$$
Then the formula (\ref{40}) implies the estimate \be\label{421}
\sup_{x\in\Re^m}f_n^\delta(x)\leq {1\over
n}\sup_{x\in\Re^m}F_n(x)+\sup_{x\in\Re^m}\int_{\Re^m}K_{n,\delta}(\|y-x\|)\mu_n(dy),
\ee where
$$
K_{n,\delta}(z)={C_1\over n}\sum_{k=1}^{[n\delta]}\brakes{n\over
k}^{m\over 2}\cdot\left[\exp\brakes{-\alpha  z^2\cdot {n\over
k}}+\brakes{1+\brakes{z\cdot \sqrt{n\over k}}^m}^{-1}\right],\quad
z\in \ax.
$$
For any $k\in \NN$ and any  $t\in\Big[{k\over n}, {k+1\over
n}\Big]$, the following inequalities  hold true:
$$
\brakes{n\over k}^{m\over 2}\leq \brakes{k+1\over k}^{m\over
2}\cdot {t}^{-{m\over 2}}\leq 2^{m\over 2}\cdot{t}^{-{m\over 2}},
$$
$$
\exp\brakes{-\alpha z^2\cdot {n\over k}}+\brakes{1+\brakes{z\cdot
\sqrt{n\over k}}^m}^{-1}\leq \exp\brakes{-{\alpha  z^2\over t}
}+\brakes{1+\brakes{{z\over \sqrt{t}}}^m}^{-1},
$$
and therefore, for any $k\in \NN$,
$$
{1\over n}\brakes{n\over k}^{m\over
2}\cdot\left[\exp\brakes{-\alpha z^2\cdot {n\over
k}}+\brakes{1+\brakes{z\cdot \sqrt{n\over k}}^m}^{-1}\right]\leq
$$
\be\label{43} \leq 2^{m\over 2}\cdot \int_{k\over n}^{k+1\over n}
{t}^{-{m\over 2}}\cdot\left[\exp\brakes{-{\alpha  z^2\over t}
}+\brakes{1+\brakes{{z\over \sqrt{t}}}^m}^{-1}\right]\, dt. \ee
The estimate (\ref{43}) implies the following estimate for the
kernel $K_{n,\delta}$:
$$
K_{n,\delta}(z)\leq 2^{m\over 2}C_1\cdot \int_{1\over
n}^\delta{t}^{-{m\over 2}}\cdot\left[\exp\brakes{-{\alpha z^2\over
t} }+\brakes{1+\brakes{{z\over \sqrt{t}}}^m}^{-1}\right]\, dt\leq
$$
$$
\leq K_\delta(z)\equiv 2^{m\over 2}C_1\cdot \int_{0}^\delta
{t}^{-{m\over 2}}\cdot\left[\exp\brakes{-{\alpha z^2\over t}
}+\brakes{1+\brakes{{z\over \sqrt{t}}}^m}^{-1}\right]\, dt
$$
We estimate the kernel $K_\delta$ using the calculations,
analogous to those made in \cite{dynkin}, \S 6.1. Making the
change of variables $u={z^2\over t}$, we get the formula
\be\label{44} K_\delta(z)=2^{m\over 2}C_1\cdot z^{2-m}\cdot
G\left({z^2\over \delta}\right), \quad G(y)=\int_y^\infty
u^{{m\over 2}-2}\left[\exp(-\alpha u)+\brakes{1+u^{m\over
2}}^{-1}\right]\,du. \ee The following properties of the function
$G$ can be verified straightforwardly: \be\label{441} G(y)\to 0,
\quad y\to +\infty,\quad
\begin{cases}G(\cdot)\hbox{ is bounded,}& m>2\\
G(y)\sim-2\ln y,\, y\to 0+,& m=2
\end{cases}.
\ee Let us give two estimates. First, \be\label{45}
K_\delta(z)\leq D\cdot w(z)\cdot G(\delta^{-{1\over 2}}),\quad
z>\delta^{1\over 4} \ee (here and below we denote by $D$ a
positive constant such that its explicit value in not needed for
us; the concrete values of $D$ may vary from line to line). For
$m>2$, this estimate follows immediately form  the representation
(\ref{44}); for $m=2$, one should recall additionally that (for
$m=2$) $w(z)\geq 1$. Second,  \be\label{46} K_\delta(z)\leq D
w(z),\quad z\in \ax.\ee For $m=3$ this estimate holds true since
the function $G$ is bounded, for $m=2$ it holds true due to the
estimates
$$
G(y)\leq D\cdot w(y),\quad  G\left({z^2\over \delta} \right)\leq
D\cdot w\left({z^2\over \delta} \right)\leq D\cdot w(z^2)\leq
2D\cdot w(z)
$$
(we used here that $\delta<1$). The formula (\ref{421}), condition
A4) and estimates (\ref{45}),(\ref{46}) imply
$$
\mathop{\lim\sup}_{n\to+\infty}\sup_{x\in\Re^m} f_n^\delta(x)\leq
D\cdot G(\delta^{-{1\over 2}})\cdot
\mathop{\lim\sup}_{n\to+\infty}\sup_{x\in\Re^m}\int_{\Re^m}w(\|y-x\|)\mu_n(dy)+
$$
\be\label{47}+D\cdot
\mathop{\lim\sup}_{n\to+\infty}\sup_{x\in\Re^m}\int_{\|y-x\|<\delta^{1\over
4}}w(\|y-x\|)\mu_n(dy). \ee Both summands in the right hand side
of the inequality (\ref{47}) tend to zero as $\delta\downarrow 0$.
This follows from (\ref{441}), Remark \ref{r22} and condition A6).
Thus, the relation (\ref{42}) is proved. The proof of the relation
(\ref{41}) is completely analogous and we omit it here.

We have verified that the processes $X_n,X$ and the functionals
$\phi_n$, defined by (\ref{00}) -- (\ref{03}), under conditions
A1) -- A6), satisfy all the conditions of Theorem \ref{t31}. Using
this theorem, we obtain the statement of  Theorem \ref{t21}.

\emph{The proof of Proposition \ref{p01}.} The family
$\{f^\delta\}$ is monotonous w.r.t. $\delta$. Therefore, the
implication (i)$\Rightarrow$(ii) follows form the Dini's Theorem
and the fact that, due to compactness of the support of $\mu$ and
the property  B2), for any given  $t\in\ax$,
$$
f^t(x)\to 0, \quad \|x\|\to+\infty.
$$
The uniform limit of a uniformly continuous functions is a
uniformly continuous function; this, together with Proposition
\ref{p41}, provides the implication (ii)$\Rightarrow$(i). The
implication (iii)$\Rightarrow$(ii) is contained in the relation
(\ref{41}). In order to prove the inverse implication, we use the
property  B3) and calculations, analogous to those made above. We
write
$$
f^\delta(x)\geq \int_{\Re^m}R_\delta(\|y-x\|)\mu(dy),\quad
R_\delta(z)\equiv\int_{0}^\delta \left[M_1{t}^{-{m\over
2}}\exp\brakes{-{\alpha z^2\over t} }-M_2{t}^{-{m\over
2}+\lambda}\exp\brakes{-{\beta z^2\over t} }\right]\, dt.
$$
Making the change of the variables  $u={z^2\over t}$, we get
$$
R_\delta(z)=z^{2-m}G_1\brakes{z^2\over
\delta}-z^{2-m+2\lambda}G_2\brakes{z^2\over \delta},
$$
$$
G_1(y)\equiv M_1\int_y^\infty u^{{m\over 2}-2}\exp[-\alpha u]\,
du,\quad G_2(y)\equiv M_2\int_y^\infty u^{{m\over
2}-2-\lambda}\exp[-\beta u]\, du.
$$
Without losing generality, we can suppose that $\lambda<1$. Let us
consider two cases. If $m\geq 3$, then the functions $G_1,G_2$ are
bounded, and, for $z<\delta$,
$$
R_\delta(z)\geq z^{2-m}\cdot
[G_1(\delta)-\delta^{2\lambda}G_2(\delta)]\geq w(z)\cdot D
$$
with $\delta$ small enough. Thus, for small enough $\delta$,
\be\label{411} f^\delta(x)\geq D \int_{\|y-x\|\leq
\delta}w(\|y-x\|)\, \mu(dy). \ee If  $m=2$, then $G_1(y)\sim -\ln
y, y\downarrow 0$, and therefore, for small enough $\delta$ and
$z<\delta$,
$$
G_1\brakes{z^2\over \delta}\geq D\ln{\delta\over z^2}\geq
D\ln{1\over z}=D w(z).
$$
Next, $G_2(y)\sim {1\over \lambda} y^{-\lambda}, y\downarrow 0$,
and thus, for small enough $\delta$ and  $z<\delta$,
$$
z^{2\lambda}G_2\brakes{z^2\over \delta}\leq
Dz^{2\lambda}\brakes{\delta\over z^2}^\lambda=D\delta^\lambda.
$$
The two latter  estimates and the fact that  $w(\cdot)\geq 1$ for
$m=2$ provide that the estimate (\ref{411}) holds true in the case
$m=2$, also. This estimate proves the implication
(ii)$\Rightarrow$(iii). The proposition is proved.

\section{Examples}

In order to illustrate the domain of application of Theorem
\ref{t21}, we give two examples. In our first example, we
construct the measure  $\mu$, satisfying condition (\ref{010}),
but not satisfying (\ref{08}). The $W$-functional, that
corresponds to this measure, can not occur as a limiting one in
the framework of the paper \cite{bass_khoshn92}; on the other
hand, this functional belongs to the class of the limiting
functionals, obtained via Theorem  \ref{t21} (see Remark
\ref{r23}).

\begin{ex}\label{e51} Let $m=2$, $\{r_k, k\geq 1\}$ be a sequence
of positive numbers (it will be defined precisely later on), and
the measure $\mu$ to have the form $\mu=\sum_{k\geq 1}Q_k\cdot
\mu_k$, where $\{Q_k, k\geq 1\}$ is a certain weight sequence and
$\mu_k=\delta_{S_k}$ is the surface measure on the circle
$S_k\equiv \{y|\|y\|=r_k\}, k\geq 1$ (the measures $\mu_k$ are
normalized in such a way that $\mu_k(S_k)=1,k\geq 1$).

For any $x\in\Re^2, \delta\in(0,e^{-1})$, denote
$B(x,\delta)=\{y|\|y-x\|\leq \delta\}$, then
$$\int_{\|y-x\|\leq
\delta}w(\|y-x\|)\mu(dy)=\sum_k Q_k\cdot V_k^\delta(x),\quad
V_k^\delta(\cdot)\equiv \int_{S_k\cap B(x,\delta)}
w(\|y-\cdot\|)\, \mu_k(dy).
$$
Note that, up to a term ${1\over 2\pi}$, $V_k^\delta$ is the
simple layer potential, generated by the measure
$\1_{B(x,\delta)}d\mu_k$, concentrated on  $S_k$
 (see, for instance, \cite{vladimirov}, \S
21). This potential is dominated by the potential $V_k$ generated
by the measure $\mu_k$, and both these potentials are a continuous
functions, harmonic in both $\{y|\|y\|<r_k\}$ and
$\{y|\|y\|>r_k\}$. The maximum principle provides that
\be\label{51}V_k(z)=V_k(0)=\ln\brakes{1\over r_k}, \, \|z\|\leq
r_k, \quad V_k^\delta(x)\leq V_k(x)\leq \ln\brakes{1\over r_k},\,
x\in\Re^2,\delta\in(0,e^{-1}),\quad k\geq 1.  \ee Due to the
maximum principle, $V_k^\delta$ takes its maximum value on the
circle $S_k$. If $\delta<r_k$, then  $S_k\cap B(x,\delta)$ is the
arch on the circle $S_k$, and it is easy to verify that the
corresponding maximal value is taken at the middle of this arch.
This reasoning and the  straightforward calculation, that is easy
and omitted, give the estimate  \be\label{52} V_k^{\delta}(x)\leq
D\cdot \brakes{\delta\over r_k}\cdot \ln\brakes {1\over \delta}.
\ee Now, we put  $r_k=2^{-k^2}, Q_k=k^{-4}, k\geq 1$. Let
$\delta\in(0,e^{-1})$ be fixed and $x\in\Re^m$ be arbitrary, let
us estimate  $\int_{\|y-x\|\leq \delta}w(\|y-x\|)\mu(dy)$. For
$\|x\|\leq 2\delta$, we have $B(x,\delta)\subset B(0,3\delta)$, an
therefore, due to  (\ref{51}), \be\label{53} \int_{\|y-x\|\leq
\delta}w(\|y-x\|)\mu(dy)\leq\sum_{k:r_k\leq 3\delta}Q_k
\ln\brakes{1\over r_k}=\sum_{k:r_k\leq 3\delta}{1\over k^2}\leq
D\cdot\left[\ln\brakes{1\over \delta}\right]^{-{1\over 2}}. \ee If
 $\|x\|>2\delta$, then there exists at most one value $k=k_*$
such that $B(x,\delta)\cap S_{k_*}\not=\emptyset$, it being known
that $r_{k_*}\geq \delta$. Then, due to  (\ref{51}),(\ref{52}),
\be\label{54} \int_{\|y-x\|\leq \delta}w(\|y-x\|)\mu(dy)\leq
D\cdot \sup_{k: r_k\geq \delta} [\ln r_k]^{-2}\cdot
\brakes{\delta\over r_k}\cdot \ln\brakes {1\over \delta}= D\cdot
\sup_{k: r_k\geq \delta}{\Phi(\delta)\over \Phi(r_k)}\cdot
\left[\ln\brakes{1\over \delta}\right]^{-1},  \ee where the
function $\Phi(t)=t[\ln t]^2, t\geq 0$ is non-decreasing in some
neighborhood of zero. From (\ref{53}),(\ref{54}), the estimate
$$
\sup_{x\in\Re^m}\int_{\|y-x\|\leq \delta}w(\|y-x\|)\mu(dy)\leq
D\cdot\left[\ln\brakes{1\over \delta}\right]^{-{1\over 2}},
$$
follows, that provides (\ref{010}). On the other hand,
$$
\mu(B(0,r_k))>Q_k=\brakes{1\over \ln{ r_k}}^2,
$$
and thus, for every $\gamma>0$,
$$
\mathop{\lim\sup}_{r\downarrow 0}r^{-\gamma}\mu(B(0,r))=+\infty,
$$
that means that condition (\ref{08}) does not hold true.
\end{ex}

In the second example, we construct a $W$-measure $\mu$, that does
not satisfy condition (\ref{010}); the $W$-functional,
corresponding to this measure, can not be obtained in Theorem
\ref{t21} as a limiting one.

\begin{ex}\label{e52} Let $m=2$, we put  $r_k=2^{-k^2},
a_k=({1\over k},0)\in\Re^m$, $\tilde S_k=\{y|\|y-a_k\|=r_k\},
\tilde Q_k=k^{-2}, k\geq 1$, $\tilde \mu_k=\delta_{\tilde S_k}$ is
the surface measure on the circle $\tilde S_k$, normalized in such
a way that  $\tilde \mu_k(\tilde S_k)=1,k\geq 1$. We put
$$
\tilde \mu=\sum_{k}\tilde Q_k\tilde \mu_k,
$$
and show that $\tilde \mu$ is a  $W$-measure that does not satisfy
(\ref{010}). By the construction, there exists $N\in \NN$ such
that, for any  $x\in\Re^m$,  the relation  \be\label{55}
\|x-a_k\|\leq 2^{-\sqrt{k}}\ee holds true for at most $N$ values
$k\in \NN$. Furthermore, if  $k>4$ and (\ref{55}) holds true, then
$\|y-x\|\leq e^{-1}, y\in \tilde S_k.$ Therefore, for a given
$x\in\Re^m$, for such $k>4$ that (\ref{55}) holds true, using the
maximum principle  we obtain the estimate
$$
\tilde Q_k\cdot \int_{\|y-x\|\leq 1}w(\|y-x\|)\tilde \mu_k(dy)\leq
\tilde Q_k\cdot\int_{\tilde S_k}\ln\brakes{1\over
\|y-a_k\|}\,\tilde \mu_k(dy)=\tilde Q_k\cdot\ln\brakes{1\over
r_k}=\tilde Q_k\cdot k^2=1.
$$
If  $k>1$ and (\ref{55}) fails, then $\|y-x\|\geq 2^{-\sqrt{k}-1},
y\in \tilde S_k.$ Thus, for a given $x\in\Re^m$, for $k>2$ such
that  (\ref{55}) fails, we have the estimate
$$
\tilde Q_k\cdot\int_{\|y-x\|\leq 1}w(\|y-x\|)\tilde \mu_k(dy)\leq
\tilde Q_k\cdot(\sqrt{k}+1)=k^{-{3\over 2}}+k^{-2}.
$$
Furthermore, every measure $\tilde\mu_k$ is a  $W$-measure, and
thus
$$
\sigma_k\equiv \sup_{x\in\Re^m}\int_{\|y-x\|\leq
1}w(\|y-x\|)\tilde \mu_k(dy)<+\infty, \quad k\geq 1.
$$
The three latter estimates imply that $$
\sup_{x\in\Re^m}\int_{\|y-x\|\leq 1}w(\|y-x\|)\tilde \mu(dy)\leq
\sum_{k=1}^4 k^{-2}\sigma_k+N+\sum_{k\geq 5}(k^{-{3\over
2}}+k^{-2})<+\infty,
$$
that means that  $\tilde \mu$ is a  $W$-measure. On the other
hand, for $\delta_k=r_k, x_k=a_k$ we have that
$$
\int_{\|y-x_k\|\leq \delta_k}w(\|y-x_k\|)\tilde \mu(dy)\geq\tilde
Q_k\cdot\int_{\|y-x_k\|\leq \delta_k}w(\|y-x_k\|)\tilde
\mu_k(dy)=1
$$
and  $\delta_k\to 0, k\to\infty.$ This means that, for $\tilde
\mu$, condition (\ref{010}) fails.
\end{ex}

\end{document}